# Generating function for Naturalized Series: The case of Ordered Motzkin Words

Gennady Eremin

ergenns@gmail.com

February 19, 2020

**Abstract.** We continue to consider the ordered lexicographic sequence, which is constructed according to the formal characteristics of a series of natural numbers. For analysis, we selected balanced parentheses with zeros, Motzkin words. As you know, generating functions allow you to work with combinatorial objects by analytical methods. Motzkin words are enumerated by Motzkin numbers, for the generation of which there is a corresponding generating function. In our case, restrictions are imposed on Motzkin words, for example, there are no leading zeros in bracket sets. The purpose of this article is to obtain the generating function of such modified Motzkin words.

A lexicographic series of ordered Motzkin words was introduced in [Er19a, Er19b]. We continue to analyze this sequence. This paper describes a generating function that generates Motzkin difference numbers – ranges of a naturalized series.

## 1   Introduction

Interest in balanced parentheses is currently quite high. Note the articles [BP14, Fan19, GZ14], in which the authors consider the partial order in parentheses. For example, in [BP14] a specific distance is established between Motzkin words in the Tamari lattice.

This article considers a *naturalized* non-numeric sequence, i.e. a sequence organized by formal attributes of a series of natural numbers. The total order allows us to introduce arithmetic, logical, and other operations for Motzkin words, up to derivatives and differential equations. As a result, we get a certain discrete analogue of classical mathematical analysis, the rudiments of *discrete calculus*.

The Motzkin word is assembled from three characters of the alphabet: zero, left (opening) parenthesis and right (closing) parenthesis. The Motzkin word can only include zeros. We will work with balanced parentheses in which a) the number of left and right parentheses is the same, and b) the number of left parentheses is not less than the number of right parentheses in every initial subword. The group of consecutive Motzkin words is the Motzkin word.

Balanced parentheses of length $n$, $n$-words, are enumerated by Motzkin numbers, $M_n$. A sequence of Motzkin numbers A001006 is known [Slo20], in which the $n$th element, $n \geq 0$, is equal to the number of Motzkin $n$-words. Here is the beginning of the sequence A001006:

(1.1)   1, 1, 2, 4, 9, 21, 51, 127, 323, 835, 2188, 5798, 15511, 41835, ...

In (1.1), elements are indexed from zero, $M_0 = 1$. Thus, a virtual *empty word* $\lambda$ of length 0, 0-word, is allowed. We know the recurrent relation for Motzkin numbers [Wei19]:

(1.2) $$M_0 = 1, \ M_n = M_{n-1} + \sum_{k=0}^{n-2} M_k M_{n-2-k}, \ n \geq 1.$$

The Motzkin word of length $n$ can start with either a left bracket or a zero (*leading zero*). Let's formulate a simple statement.

**Proposition 1.1.** *The number of Motzkin n-words, $n \geq 2$, that begin with the leading zero, is $M_{n-1}$.*

*Proof.* If we add a leading zero to each Motzkin word of length $n-1$, $n \geq 2$, we get $M_{n-1}$ Motzkin $n$-words. We can add a left (or right) parenthesis in front of the Motzkin $(n-1)$-word, but this will not be the Motzkin word due to the unbalance of parentheses. □

Let's call balanced parentheses with leading zero *inherited*. The Motzkin words starting with the left bracket are called *unique*. The number of unique Motzkin $n$-words will be denoted by $U_n$. We will also consider unique a single word of length 1, the Motzkin word "0", so $U_1 = 1$. An empty word cannot be unique, therefore $U_0 = 0$. Based on (1.2), we can write the Motzkin *difference numbers* as follows [Ber99]:

(1.3) $$U_0 = 0, \ U_1 = 1, \ U_n = M_n - M_{n-1} = \sum_{k=0}^{n-2} M_k M_{n-2-k}, \ n \geq 2.$$

Many mathematicians consider zero to be a natural number. This suits us, because among the words there is an identical Motzkin 1-word. Based on the formal properties of natural numbers, as well as the order in the alphabet $0 < ( < )$, we get the *naturalized series* of unique Motzkin words [Er19b]:

$\mathfrak{M} = \{0,\ (),\ (0),\ ()0,\ (00),\ (0)0,\ (()),\ ()00,\ ()(),\ (000),\ (00)0,\ (0()),\ ...\}$

Elements are indexed from zero and grouped by code length (initial words in ranges are highlighted in red): $\mathfrak{M}_1 = \{0\}$, $_2 = \{()\}$, $\mathfrak{M}_3 = \{(0), ()0\}$, and so on. The power ranges are as follows: $\#\mathfrak{M}_1 = U_1 = 1$, $\#\mathfrak{M}_n = U_n = M_n - M_{n-1}$, $n > 1$. We work with the real Motzkin words, therefore $\#\mathfrak{M}_0 = U_0 = 0$. The difference numbers $U_n$ for $n = 0, 1, 2, ...$ form the sequence:

(1.4)   0, 1, 1, 2, 5, 12, 30, 76, 196, 512, 1353, 3610, 9713, 26324, 71799, ….

In the next section, we calculate the generating function of the sequence (1.4), in other words, a formal power series whose coefficients are the difference numbers $U_n$.

## 2   The generating function for the Motzkin difference numbers

Generating functions allow you to work with combinatorial objects by analytical methods. Often generating functions help to derive explicit formulas for the number of combinatorial objects [FS09, LM12, Lan03].

**2.1.** In recent years, the symbolic method of constructing generating functions has been popular. The feature of the symbolic method is the use of inference rules for the language grammar.

We will use a two-step procedure. Let's start with the rules for deriving the Motzkin language. The Motzkin word is either

– the empty word λ of length 0, or
– the word $(a)\,b$, where $a$ and $b$ are the Motzkin words, or
– the word $0a$, where $a$ is the Motzkin word.

The first two rules coincide the output rules of the Dyck language. Let $\mathcal{M}$ be the set of all Motzkin words. Three rules correspond to the structural equation

$$\mathcal{M} = \lambda + (\mathcal{M})\mathcal{M} + 0\mathcal{M},$$

where plus denotes the union of disjoint sets. Further, replacing the set $\mathcal{M}$ by the generating function $M(x)$, we get the functional equation [DS77]:

(2.1)   $M(x) = 1 + x^2 M^2(x) + xM(x)$

$= \sum_{n \geq 0} M_n x^n = 1 + x + 2x^2 + 4x^3 + 9x^4 + 21x^5 + 51x^6 + \ldots$

Note that in (2.1) the alphabet characters (, 0, and ) are replaced by $x$, and the empty word λ corresponds to 1. An empty word is the only one, so the number of Motzkin words increases by exactly 1.

When decomposing the generating function into a formal Taylor power series, we get Motzkin numbers, in other words, these numbers are generated by $M(x)$. The solution for (2.1) is as follows [Wei19]:

(2.2)   $$M(x) = \frac{1-x-\sqrt{1-2x-3x^2}}{2x^2} = \frac{2}{1-x+\sqrt{1-2x-3x^2}};$$

The latter expression is more convenient, because division by zero is eliminated (at the point $x = 0$).

Let's move on to the second stage, to our naturalized series 𝔐. We need to get a generating function for the difference numbers $U_0 = 0$, $U_1 = 1$, and $U_n = M_n - M_{n-1}$, $n > 1$. Recall that an empty word is not a unique Motzkin word, so $U_0 = 0$. And this is logical, because we are building a sequence as close as possible to natural numbers, and there is also no empty word of length 0 among natural numbers (for example, the smallest integer 0 is real and has a length of 1).

In series 𝔐, there are no Motzkin words with leading (initial) zeros. The only 1-word "0" starts with zero. This single case resembles an empty Motzkin word (see the first rule in the Motzkin language), so it is logical to formulate a similar rule for unique Motzkin words. Other unique words begin with the left parenthesis, so it's enough to repeat the second rule of Motzkin language.

Thus, in the second stage, we formulate two additional rules: the unique Motzkin word is either

– the word "0", or
– the word $(a)b$, where $a$ and $b$ are ordinary Motzkin words.

As a result, we obtain a structural equation for the set of unique Motzkin words

$$\mathfrak{M} = \text{'0'} + (M)M.$$

Then the generating function for the Motzkin difference numbers takes the form

(2.3) $$\text{Nat}(x) = x + x^2 M^2(x),$$

where $M(x)$ is the generating function of ordinary Motzkin words. We again replaced the alphabet symbols (, 0, and ) with $x$.

**2.2.** The required generating function (2.3) has been obtained, and we could proceed to the calculations. But as it often happens, there is a desire to check at first glance an unusual and somewhat strange symbolic method of constructing generating functions. In this case, it is not difficult to get a duplicate of (2.3) in the traditional way. Let's use the formula (1.3)

$$\text{Nat}(x) = \sum_{n \geq 0} U_n x^n = U_0 + U_1 x + \sum_{n \geq 2} (M_n - M_{n-1}) x^n$$
$$= x + \sum_{n \geq 2} M_n x^n - \sum_{n \geq 2} M_{n-1} x^n;$$

$$x + \sum_{n \geq 2} M_n x^n = x + M_0 + M_1 x + \sum_{n \geq 2} M_n x^n - M_0 - M_1 x$$
$$= x + \sum_{n \geq 0} M_n x^n - 1 - x = M(x) - 1;$$

$$\sum_{n \geq 2} M_{n-1} x^n = x \sum_{n \geq 2} M_{n-1} x^{n-1} = x(M_1 x + M_2 x^2 + \ldots + M_0 - M_0)$$
$$= x(M(x) - 1).$$

As a result, we get

(2.4) $\quad \text{Nat}(x) = M(x) - 1 - x(M(x) - 1) = x - 1 + (1 - x)M(x).$

At first glance, the formulas (2.3) and (2.4) differ, but this is not so. It is enough to go back to (2.1) and get the obvious identity

$$x - 1 + M(x) - xM(x) = x + x^2 M^2(x).$$

**2.3.** In future calculations, we will use both formulas (2.3) and (2.4). After substitution (2.2) we get

(2.5) $\quad \text{Nat}(x) = x + x^2 \left(\dfrac{2}{1-x+\sqrt{1-2x-3x^2}}\right)^2 = x - 1 + \dfrac{2 - 2x}{1-x+\sqrt{1-2x-3x^2}}$

$\qquad = \sum_{n \geq 0} U_n x^n = 0 + x + x^2 + 2x^3 + 5x^4 + 12x^5 + 30x^6 + 76x^6 + \ldots$

From the first sum (two terms with squares), we obtain the initial coefficients

$$U_0 = \text{Nat}(0)/0! = 0/1 = 0, \quad U_1 = \text{Nat}'(0)/1! = 1/1 = 1.$$

The second option (the fraction without squares, three terms) is convenient for further calculations. Here it is enough to consider only the last term (the beginning of $x - 1$ is reset after the second derivative). Let's consider a small algorithm.

**Algorithm 2.1.** Let's denote: $s = 1 - 2x - 3x^2$, $t = s'/2 = -1 - 3x$, $\sqrt{} = \sqrt{s}$, $\sqrt{}\sqrt{} = s$, $\sqrt{}' = s'/2\sqrt{} = (-1-3x)/\sqrt{} = t/\sqrt{}$. Obviously, $\sqrt{}\sqrt{}' = t$. Write the fraction as

$$U/V = (2-2x)/(1 - x + \sqrt{}) = (a + b\sqrt{})/(c + d\sqrt{}),$$

with initial polynomial values: $a = 2 - 2x$, $b = 0$, $c = 1 - x$, $d = 1$.

Then we calculate the first derivative $(U/V)' = (U'V - UV')/V^2$:

$U'V = V(a + b\sqrt{})' = Va' + Vb'\sqrt{} + Vb\sqrt{}' = (c + d\sqrt{})a' + (c + d\sqrt{})b'\sqrt{} + (c + d\sqrt{})b\sqrt{}'$

$\qquad = (a'c + b'ds + bdt) + (a'd + b'c)\sqrt{} + bct\sqrt{}^{-1}$;

$UV' = U(c + d\sqrt{})' = Uc' + Ud'\sqrt{} + Ud\sqrt{}' = (a + b\sqrt{})c' + (a + b\sqrt{})d'\sqrt{} + (a + b\sqrt{})d\sqrt{}'$

$\qquad = (ac' + bd's + bdt) + (bc' + ad')\sqrt{} + adt\sqrt{}^{-1}$;

$V^2 = c^2 + d^2\sqrt{}\sqrt{} + 2cd\sqrt{} = c^2 + d^2 s + 2cd\sqrt{}$.

As a result, we get a new fraction $(U/V)' = (A + B\sqrt{})/(C + D\sqrt{})$, where

$A = (a'd + b'c - bc' - ad')s + (bc - ad)t$; $B = a'c - ac' + (b'd - bd')s$;
$C = 2cds$; $D = c^2 + d^2 s$.

It remains for us to organize a cycle, i.e. repeat the calculations after updating the polynomials. In General, on the $k$-th pass of the cycle, $k > 1$, we get the coefficient:

$U_k = \text{Nat}^{(k)}(0)/k! = (A|_{x=0} + B|_{x=0})/((C|_{x=0} + D|_{x=0})k!).\qquad \square$

The JScript program is available to the reader; you can check the calculations using the link https://eremin.000webhostapp.com/arxiv/gener-fun.html. The usual bit grid of a computer allows you to get up to $U_{12} = 9713$.

**Acknowledgements.** I would like to thank Sergey Kirgizov (LIB, Univ. Bourgogne Franche-Comté, France) for his help in calculating the generating functions.

Gzhel State University, Moscow, 140155, Russia
http://www.en.art-gzhel.ru/